\documentclass[preprint]{elsarticle}
\usepackage[symbol]{footmisc}

\usepackage{mathrsfs}
\usepackage{amssymb}
\usepackage{amsmath}
\usepackage{graphics}
\usepackage{latexsym}
\usepackage{amsfonts}
\usepackage{color}
\usepackage{epstopdf}
\usepackage{ctable}
\usepackage{graphicx}
\usepackage{float} 

\usepackage{verbatim}
\newtheorem{theo}{Theorem}[section]

\newdefinition{example}[theo]{Example}
\newtheorem{uppg}[theo]{Exercise}
\newdefinition{remark}[theo]{Remark}
\newdefinition{definition}[theo]{Definition}
\newproof{pf}{Proof}
\newproof{pot}{Proof of Theorem \ref{thm2}}
\newcommand{\be}{\begin{eqnarray*}}
\newcommand{\ee}{\end{eqnarray*}}
\newcommand{\ben}{\begin{eqnarray}}
\newcommand{\een}{\end{eqnarray}}

\def\subsecn (#1) {\medskip\ \ \ {\it #1}\medskip}

\usepackage{mathrsfs}
\usepackage{amssymb}
\usepackage{array}
\usepackage{multirow}
\usepackage{makecell}
\newcommand{\lp}[1]{\left(\begin{array}{#1}}
\newcommand{\rp}{\end{array}\right)}

\newcommand{\leftd}[1]{\left\{\begin{array}{#1}}
\newcommand{\rightd}{\end{array}\right.}

\def\E {\mathbf{E}}

\date{\today}

\begin{document}

\begin{frontmatter}
\title{Bayesian estimation of a competing risk model based on Weibull and exponential distributions under right censored data} 

\author[talhi]{Hamida Talhi}
\author[hib]{Hiba Aiachi}
\author[nadji]{Nadji Rahmania\fnref{fn1}}
\address{$^{1}$Probability Statistics Laboratory, Badji Mokhtar University, BP12, 23000 Annaba, Algeria }
\address{$^{1,1}$Paul Painlev\'e laboratory, UMR-CNRS 8524.
Lille University, 59655 Villeneuve d'Ascq C\'edex, France.}
\fntext[fn1]{Corresponding Author.\\E-mail adress: nadji.rahmania@univ-lille.fr (N.Rahmania).}
\begin{abstract}
In this paper we investigate the estimation of the unknown parameters of a competing risk model based on a Weibull distributed decreasing failure rate and an exponentially distributed constant failure rate, under right censored data. 
The Bayes estimators and the corresponding risks are derived using various loss functions. Since the posterior analysis involves analytically intractable integrals, we propose a Monte-Carlo method to compute these estimators. Given initial values of the model parameters, the Maximum Likelihood estimators are computed using  the Expectation-Maximization algorithm. Finally, we use Pitman's closeness criterion and integrated mean-square error to compare the performance of the Bayesian and the maximum likelihood estimators.
\end{abstract}


\end{frontmatter}


\section{Introduction}

The exponential and the Weibull distributions are the most used distributions in life time data analysis, mostly due to experience and goodness-of-fit tests,  see for instance Lawless (2002) and Hamada et al. (2008). In this paper we propose a Bayesian analysis of a computing risk model based on Weibull and exponential distributions under right censored data. Boudjerda et al. (2016) considered the Bayesian analysis of the right truncated Weibull distribution under type II censored data and derived Bayes estimators and the corresponding risks using symmetric and asymmetric loss functions. Aouf and Chadli (2017) considered the Bayesian analysis of generalized Lindley distribution under type II censored data and derived Bayes estimators and the corresponding risks using symmetric and asymmetric loss functions. Balakrishnan and Mitra (2012) applied the EM algorithm to estimate the parameters of the Weibull distribution when the model is left-truncated and the data are right censored.

\noindent The exponential distribution $\mathcal{E}(\eta)$, with mean $\eta$ is often used for modelling failure times caused by accidents cleared of birth defects of a no ageing material. The survival function of the exponential distribution is
\[S_E(t)=\exp (-\frac{t}{\eta}),\]
where the scale parameter $\eta$ is the inverse of the constant hazard rate $\lambda$. 

\noindent The versatile Weibull $\mathcal{W}(\eta,\beta)$ distribution, has survival function
\[S_W(t)=\exp[-\left(\frac{t}{\eta}\right)^\beta]\]
and hazard rate
\[h_W(t)=\frac{\beta}{\eta}\left(\frac{t}{\eta}\right)^{\beta-1}.\]
When the shape parameter $\beta<1$, the decreasing hazard rate of the model can be used for modelling failure due to early birth defects, and when $\beta>1$ in can be used as a model failures due to ageing. When $\beta=1$, the Weibull distribution reduces to the exponential distribution with scale parameter $\eta$. This last case may arise due to failures by accidents. 

When modeling reliability feedback data with the Weibull distribution, the problem one is faced with is  to decide whether $\beta=1$ or $\beta<1$ when the failure is due to birth defects or $\beta=1$ versus $\beta>1$ when the failure is due to ageing. This question can be solved by using likelihood ratio tests as suggested in e.g. d'Agostino and Stephens (1986). For simplicity, we only consider failures due to ageing only. If, for instance, ageing is diagnosed, then further statistical inference is made assuming that the observed failure times arise from a Weibull distribution.  This assumption may be plausible in many circumstances, but there are many situations where neglecting accidental failures might introduce an non-negligible bias in statistical inference on material lifetimes. Even when ageing is the most frequent cause of failure, accidental failures may still  be numerous. Thus, a more realistic way of modelling failure times is a competing risk model which takes into account the fact that a failure can be caused by ageing or accidents. This later model is defined as follows.
A failure time is the realisation of the random variable $B=min(E,W)$ where $E$ follows the  exponential distribution $\mathcal{E}(\eta_0)$  and $W$ follows the Weibull distribution $W(\eta_1,\beta)$ where $\beta>1$, where the r.v.  $E$ and $W$ are assmed to be independent. Consequently, the distribution of $B$ is characterised by the  parameters $\eta_0,\eta_1$ and $\beta$ and will be denoted $\mathcal{B}(\eta_0,\eta_1,\beta)$.

We propose two approaches to the estimate the parameters  of $\mathcal{B}(\eta_0,\eta_1,\beta)$. The first approach is the classical maximum likelihood estimation (MLE) and the second one is the Bayesian estimation using three loss functions (Generalized Quadratic function, entropy function  and Linex function). We use the Metropolis Hastings sampling procedure to generate Monte-Carlo samples to obtain the Bayes estimators of the unknown parameters. Finally, we perform some simulation experiments to compare the performance of the proposed Bayes estimators and the maximum likelihood estimators in terms of Pitman's closeness criterion and the integrated mean square error (IMSE).

The rest of the the paper is organized as the following: In Section 2, we present the main characteristics of the model. Section 3 deals with the maximum likelihood estimation of the $\mathcal{B}(\eta_0,\eta_1,\beta)$ distribution through the EM algorithm. In Section 4, the Bayesian estimators under different loss functions are displayed. Monte-Carlo simulation results are presented in Section 5. Finally, Section 6  concludes the paper. 

\section{The $\mathcal{B}$ distribution}
Consider the r.v. $B=min(E,W)$ where $E$ is exponentially distributed with mean $\eta_0$ and $W$ follows the Weibull distribution with scale parameter $\eta_1$ and shape parameter $\beta$, $E$ and $W$ being independent. The main characteristics of the probability distribution $\mathcal{B}$  of the r.v. $B$ are as follows. Its hazard function is
\begin{equation}\label{h-B}
h_B(x)=\frac{1}{\eta_0}+\frac{\beta}{\eta_1}\left(\frac{x}{\eta_1}\right)^{\beta-1},
\end{equation}
its survival (or reliability) function is
\begin{equation} \label{S-B}
S_B(x)=\text{exp}\left[-\frac{x}{\eta_0}-\left(\frac{x}{\eta_1}\right)^{\beta}\right]
\end{equation}
and its probability density function (pdf) is
\begin{equation}\label{f-B}
f_B(x)=\left(\frac{1}{\eta_0}+\frac{\beta}{\eta_1}\left(\frac{x}{\eta_1}\right)^{\beta-1}\right)\text{exp}\left[-\frac{x}{\eta_0}-\left(\frac{x}{\eta_1}\right)^{\beta}\right].
\end{equation}
\section{Maximum likelihood estimation}
Consider a $n$-sample $(X_1,X_2,\dots,X_n)$ generated from the $\mathcal{B}$ distribution with pdf \eqref{f-B}. Assuming the data is right censored,  the likelihood function for right censoring data reads
\begin{align*}
L(\eta_0,\eta_1,\beta|X)&=\prod_{i=1}^{n}f_B(x_i)^{\delta_i}S_B(x_i)^{1-\delta_i}=\prod_{i=1}^{n}h_B(x_i)^{\delta_i}S_B(x_i),\quad x_1\leq x_2\leq \cdots \leq x_n.
\end{align*}
In view of \eqref{h-B} and \eqref{S-B}, the likelihood function is
\begin{equation}
L(\eta_0,\eta_1,\beta|X)=\prod_{i=1}^n\left(\frac{1}{\eta_0}+\frac{\beta}{\eta_1}\left(\frac{x_i}{\eta_1}\right)^{\beta-1}\right)^{\delta_i}\text{exp}\left[-\frac{\sum_{i=1}^n x_i}{\eta_0}-\sum_{i=1}^n\left(\frac{x_i}{\eta_1}\right)^{\beta}\right].
\end{equation}

Since the r.v. $B$ is a result of a competition between $E$ and $W$,  the data model is incomplete in the sense that although the observations are realizations of one of these r.v.,  it is often hard to know beforehand whether a particular observation is a realization of $E$ or $W$. This in turn makes a direct maximization of the likelihood function numerically highly unstable. Instead, we use the EM algorithm with its two steps, expectation (E) and maximization (M), seems a plausible  alternative to the direct maximization of the likelihood function for incomplete data models, especially when we can implement the maximization step separately for the exponential and the Weibull models (cf. Dempster et al. (1977), Bousquet et al. (2006) and Little \& Rubin (2002)).

We proceed as follows. Define $z_i=(z^E_i,z_i^W)$ where $z^E_i=1$ and $z^W_i=0$ indicates that the associated observation is coming from an exponential model, and $z^W_i=1$ and $z^E_i=0$ is from the Weibull distribution. By convention, the complete data can then be written as $o=(o_i=(x_i,z_i), i=1...,n)=(x,z).$ So the resulting competing risk density can be written as
\begin{equation}
f(o_i)=h_E(x_i)^{z^E_i}h_W(x_i)^{z^W_i}S_E(x_i)S_W(x_i),
\end{equation}
and the log-likelihood based on complete data $o=(o_1,...,o_n)$ reads
\begin{equation}
l(\eta_0,\eta_1,\beta \, | o)=\sum_{i=1}^n z^E_i\log(h_E(x_i))+z^W_i\log(h_W(x_i))+\log(S_E(x_i))+\log(S_W(x_i)).
\end{equation}
Set $\Theta=(\eta_0,\eta_1,\beta)$ and let $\tilde{\Theta}$ denote its current value. The expected value of log-likelihood $Q(\Theta|\tilde{\Theta})$ is 
\begin{align}\label{Q}
Q(\Theta|\tilde{\Theta})&=E(l(\eta_0,\eta_1,\beta|o)|x,\tilde{\Theta})\nonumber\\
&=\sum_{i=1}^n \tilde{p}_E(x_i)\log(h_E(x_i))+\tilde{p}_W(x_i)\log(h_W(x_i))+\log(S_E(x_i))+\log(S_W(x_i)),
\end{align}
where \[\tilde{p}_E(x_i)= E(\tilde{p}_E(x_i|x,\tilde{\Theta})=P(z^E_i=1|x,\tilde{\Theta})=\frac{h_E(x_i)}{h_E(x_i)+h_W(x_i)}, \quad \tilde{p}_W(x_i)=1-\tilde{p}_E(x_i).\]
Here $\tilde{p}_E(x_i)\;(\tilde{p}_W(x_i))$ denotes the probability that the observation is coming
from the exponential (Weibull) distribution. Moreover, the equation \eqref{Q} has an additive structure that results from the contribution of both
the exponential and the Weibull distributions. This additive decomposition of \eqref{Q} makes the implementation of M-step easier in the sense that it maximizes separately the terms corresponding to exponential and Weibull distributions. The exponential term can be maximized by direct differentiation with respect to the parameter $\eta_0$, whereas the Weibull term can be maximized using any of the iterative procedures such as  the Newton-Raphson method (see Mann et al. (1974), Press et al. (2007)),  since there is no closed form of the derivatives with respect to the Weibull parameters $\eta_1$ and $\beta$. These two steps can be repeated until the iterating algorithm converges to give the desired MLE estimates.

\section{Bayesian estimators under different loss functions}
In the Bayesian approach, the unknown parameters are considered as random variables (r.v) instead of fixed constants, from this point the variations in the parameters can be incorporated by assuming prior distributions of the unknown parameters. As prior distributions,  we assume  the parameters $\eta_0, \eta_1$  follow the Gamma distribution as a prior:
\begin{align*}
\pi(\eta_0)=&\frac{a_1^{b_1}}{\Gamma(b_1)}\eta_0^{b_1-1}\text{exp}[-a_1\eta_0]\\
\pi(\eta_1)=&\frac{a_2^{b_2}}{\Gamma(b_2)}\eta_1^{b_2-1}\text{exp}[-a_2\eta_1],\\
\end{align*}
while the parameter $\beta$ follow an uniform distribution, $\beta\sim U(\beta_l,\beta_r)$,
 \[\pi(\beta)=\frac{1}{\beta_r-\beta_l}, \quad \beta_r\leq\beta\leq\beta_l.\]
Moreover, $\eta_0,\eta_1$ and $\beta$ are assumed  independent.  Thus, the joint prior distribution of $(\eta_0,\eta_1,\beta)$ is given by
\begin{equation}
\pi(\eta_0,\eta_1,\beta)=\frac{a_1^{b_1}a_2^{b_2}}{\Gamma(b_1)\Gamma(b_2)}\eta_0^{b_1-1}\eta_1^{b_2-1}\text{exp}[-a_1\eta_0-a_2\eta_1]\frac{1}{\beta_r-\beta_l}.
\end{equation}
There is no specific criterion for the selection of the Gamma family except that it is flexible and admits a Gamma distribution as a conjugate prior. The posterior density is then 
\begin{equation*}
\pi(\eta_0,\eta_1,\beta|X)=\frac{L(\eta_0\eta_1,\beta|X)\pi(\eta_0,\eta_1,\beta)}{\int\int\int_0^{+\infty}L(\eta_0\eta_1,\beta|X)\pi(\eta_0,\eta_1,\beta)\text{d}\eta_0\text{d}\eta_1\text{d}\beta},
\end{equation*}
so the joint posterior of $(\eta_0,\eta_1,\beta)$ is
\begin{eqnarray}
\pi(\eta_0,\eta_1,\beta|X)=K\eta_0^{b_1-1}\eta_1^{b_2-1}e^{\left[-\frac{\sum_{i=1}^n x_i}{\eta_0}-\sum_{i=1}^n\left(\frac{x_i}{\eta_1}\right)^{\beta}-a_1\eta_0-a_2\eta_1\right]}\prod_{i=1}^n\left(\frac{1}{\eta_0}+\frac{\beta}{\eta_1}\left(\frac{x_i}{\eta_1}\right)^{\beta-1}\right)^{\delta_i},
\end{eqnarray}
where $K$ is the normalizing constant. 

\noindent Next, we introduce the three loss functions namely  the generalised quadratic (GQ), the Linex and the entropy functions we will consider below. In  the following table we display these loss functions with their Bayes estimators and the corresponding posterior risks (PR).
\begin{table}[htp]
\footnotesize\setlength{\tabcolsep}{4pt}
\begin{tabular}{*{4}{|c}|}
\hline
Loss function & Expression & Bayes estimators & posterior risk \\
\hline
Generalised quadratic & $L(\lambda,\delta)=\tau(\lambda)(\lambda-\delta)^2$ & $\hat{\delta}_{GQ}=\frac{E_{\pi}(\tau(\lambda)\lambda)}{E_{\pi}(\tau(\lambda)}$ &$ E_{\pi}(\tau(\lambda)(\lambda-\delta)^2 $ \\
\hline
Entropy & $L(\lambda,\delta)=\left(\frac{\delta}{\lambda}\right)^p-p\log\left(\frac{\delta}{\lambda}\right)-1 $&$ \hat{\delta}_{E}=E_{\pi}(\lambda^{-p})^{\frac{-1}{p}}$ & $p[E_{\pi}(\log(\lambda-\log(\hat{\delta}_E)))]$\\
\hline
Linex & $L(\lambda,\delta)=\exp(r(\delta-\lambda))-r(\delta-\lambda)-1$ & $\hat{\delta}_{L}=\frac{-1}{r}\log(E_{\pi}(\exp(-r\lambda))$ & $r(\hat{\delta}_{GQ}-\hat{\delta}_{L})$\\
\hline
\end{tabular}
\caption{The loss functions and the corresponding Bayesian estimators and the posterior risk}
\end{table}
Under the GQ loss function, $L(\lambda,\delta)=\tau(\lambda)(\lambda-\delta)^2$,
assuming that $\tau(\lambda)=\lambda^{\alpha-1}$, the Bayesian estimators of $\eta_0$, $\eta_1$ and $\beta$ denoted respectively by $\hat{\eta}_{0(GQ)}, \hat{\eta}_{1(GQ)}$ and $\hat{\beta}_{(GQ)}$ are 
\begin{table}[H]
\footnotesize\setlength{\tabcolsep}{4pt}
\begin{align*}
\hat{\eta}_{0(GQ)}&=\frac{\int\int\int_0^{+\infty}\eta_0^{\alpha+b_1-1}\eta_1^{b_2-1}\text{exp}\left[-\frac{\sum_{i=1}^n x_i}{\eta_0}-\sum_{i=1}^n\left(\frac{x_i}{\eta_1}\right)^{\beta}-a_1\eta_0-a_2\eta_1\right]\prod_{i=1}^n\left(\frac{1}{\eta_0}+\frac{\beta}{\eta_1}\left(\frac{x_i}{\eta_1}\right)^{\beta-1}\right)^{\delta_i}\text{d}\eta_0\text{d}\eta_1\text{d}\beta}{\int\int\int_0^{+\infty}\eta_0^{\alpha+b_1-2}\eta_1^{b_2-1}\text{exp}\left[-\frac{\sum_{i=1}^n x_i}{\eta_0}-\sum_{i=1}^n\left(\frac{x_i}{\eta_1}\right)^{\beta}-a_1\eta_0-a_2\eta_1\right]\prod_{i=1}^n\left(\frac{1}{\eta_0}+\frac{\beta}{\eta_1}\left(\frac{x_i}{\eta_1}\right)^{\beta-1}\right)^{\delta_i}\text{d}\eta_0\text{d}\eta_1\text{d}\beta},\\
\hat{\eta}_{1(GQ)}&=\frac{\int\int\int_0^{+\infty}\eta_0^{b_1-1}\eta_1^{\alpha+b_2-1}\text{exp}\left[-\frac{\sum_{i=1}^n x_i}{\eta_0}-\sum_{i=1}^n\left(\frac{x_i}{\eta_1}\right)^{\beta}-a_1\eta_0-a_2\eta_1\right]\prod_{i=1}^n\left(\frac{1}{\eta_0}+\frac{\beta}{\eta_1}\left(\frac{x_i}{\eta_1}\right)^{\beta-1}\right)^{\delta_i}\text{d}\eta_0\text{d}\eta_1\text{d}\beta}{\int\int\int_0^{+\infty}\eta_0^{b_1-1}\eta_1^{\alpha+b_2-2}\text{exp}\left[-\frac{\sum_{i=1}^n x_i}{\eta_0}-\sum_{i=1}^n\left(\frac{x_i}{\eta_1}\right)^{\beta}-a_1\eta_0-a_2\eta_1\right]\prod_{i=1}^n\left(\frac{1}{\eta_0}+\frac{\beta}{\eta_1}\left(\frac{x_i}{\eta_1}\right)^{\beta-1}\right)^{\delta_i}\text{d}\eta_0\text{d}\eta_1\text{d}\beta},\\
\hat{\beta}_{(GQ)}&=\frac{\int\int\int_0^{+\infty}\beta^{\alpha}\eta_0^{b_1-1}\eta_1^{b_2-1}\text{exp}\left[-\frac{\sum_{i=1}^n x_i}{\eta_0}-\sum_{i=1}^n\left(\frac{x_i}{\eta_1}\right)^{\beta}-a_1\eta_0-a_2\eta_1\right]\prod_{i=1}^n\left(\frac{1}{\eta_0}+\frac{\beta}{\eta_1}\left(\frac{x_i}{\eta_1}\right)^{\beta-1}\right)^{\delta_i}\text{d}\eta_0\text{d}\eta_1\text{d}\beta}{\int\int\int_0^{+\infty}\beta^{\alpha-1}\eta_0^{b_1-1}\eta_1^{b_2-1}\text{exp}\left[-\frac{\sum_{i=1}^n x_i}{\eta_0}-\sum_{i=1}^n\left(\frac{x_i}{\eta_1}\right)^{\beta}-a_1\eta_0-a_2\eta_1\right]\prod_{i=1}^n\left(\frac{1}{\eta_0}+\frac{\beta}{\eta_1}\left(\frac{x_i}{\eta_1}\right)^{\beta-1}\right)^{\delta_i}\text{d}\eta_0\text{d}\eta_1\text{d}\beta}.
\end{align*}
\end{table}
\noindent The corresponding posterior risks are then
\[PR(\hat{\eta}_{0(GQ)})=E_{\pi}(\eta_0^{\alpha+1})-2\hat{\eta}_{0(GQ)}E_{\pi}(\eta_0^{\alpha})+\hat{\eta}_{0(GQ)}^2E_{\pi}(\eta_0^{\alpha-1}),\]
\[PR(\hat{\eta}_{1(GQ)})=E_{\pi}(\eta_1^{\alpha+1})-2\hat{\eta}_{1(GQ)}E_{\pi}(\eta_1^{\alpha})+\hat{\eta}_{1(GQ)}^2E_{\pi}(\eta_1^{\alpha-1}),\]
\[PR(\hat{\beta}_{(GQ)})=E_{\pi}(\beta^{\alpha+1})-2\hat{\beta}_{(GQ)}E_{\pi}(\beta^{\alpha})+\hat{\beta}_{(GQ)}^2E_{\pi}(\beta^{\alpha-1}).\]
We note that when $\alpha=1$, we retrieve the basic quadratic loss function.

\medskip Under the entropy loss function, the Bayesian estimators $\hat{\eta}_{0(E)}$, $\hat{\eta}_{0(E)}$ and $\hat{\beta}_{(E)}$ are 
\begin{table}[H]
\footnotesize\setlength{\tabcolsep}{4pt}
\begin{align*}
\hat{\eta}_{0(E)}&=\left[K\int\int\int_0^{+\infty}\eta_0^{b_1-1-p}\eta_1^{b_2-1}\text{exp}\left[-\frac{\sum_{i=1}^n x_i}{\eta_0}-\sum_{i=1}^n\left(\frac{x_i}{\eta_1}\right)^{\beta}-a_1\eta_0-a_2\eta_1\right]\prod_{i=1}^n\left(\frac{1}{\eta_0}+\frac{\beta}{\eta_1}\left(\frac{x_i}{\eta_1}\right)^{\beta-1}\right)^{\delta_i}\text{d}\eta_0\text{d}\eta_1\text{d}\beta\right]^{\frac{-1}{p}},\\
\hat{\eta}_{1(E)}&=\left[K\int\int\int_0^{+\infty}\eta_0^{b_1-1}\eta_1^{b_2-1-p}\text{exp}\left[-\frac{\sum_{i=1}^n x_i}{\eta_0}-\sum_{i=1}^n\left(\frac{x_i}{\eta_1}\right)^{\beta}-a_1\eta_0-a_2\eta_1\right]\prod_{i=1}^n\left(\frac{1}{\eta_0}+\frac{\beta}{\eta_1}\left(\frac{x_i}{\eta_1}\right)^{\beta-1}\right)^{\delta_i}\text{d}\eta_0\text{d}\eta_1\text{d}\beta\right]^{\frac{-1}{p}},\\
\hat{\beta}_{(E)}&=\left[K\int\int\int_0^{+\infty}\beta^{-p}\eta_0^{b_1-1}\eta_1^{b_2-1}\text{exp}\left[-\frac{\sum_{i=1}^n x_i}{\eta_0}-\sum_{i=1}^n\left(\frac{x_i}{\eta_1}\right)^{\beta}-a_1\eta_0-a_2\eta_1\right]\prod_{i=1}^n\left(\frac{1}{\eta_0}+\frac{\beta}{\eta_1}\left(\frac{x_i}{\eta_1}\right)^{\beta-1}\right)^{\delta_i}\text{d}\eta_0\text{d}\eta_1\text{d}\beta\right]^{\frac{-1}{p}}.
\end{align*}
\end{table}
\noindent The corresponding posterior risks are then
\[PR(\hat{\eta_0}_{(E)})=pE_{\pi}(\log(\eta_0)-\log(\hat{\eta}_{0(E)})),\]
\[PR(\hat{\eta}_{1(E)})=pE_{\pi}(\log(\eta_1)-\log(\hat{\eta}_{1(E)})),\]
\[PR(\hat{\beta}_{(E)})=pE_{\pi}(\log(\beta)-\log(\hat{\beta}_{(E)})).\]

Finally,  under the Linex loss function we obtain the following estimators
\begin{table}[H]
\scriptsize{
\begin{align*}
\hat{\eta}_{0(L)}&=\frac{-K}{r}\log\left[\int\int\int_0^{+\infty}\eta_0^{b_1-1}\eta_1^{b_2-1}e^{\left[-\frac{\sum_{i=1}^n x_i}{\eta_0}-\sum_{i=1}^n\left(\frac{x_i}{\eta_1}\right)^{\beta}-a_1\eta_0-a_2\eta_1-r\eta_0\right]}\prod_{i=1}^n\left(\frac{1}{\eta_0}+\frac{\beta}{\eta_1}\left(\frac{x_i}{\eta_1}\right)^{\beta-1}\right)^{\delta_i}\text{d}\eta_0\text{d}\eta_1\text{d}\beta\right],\\
\hat{\eta}_{1(L)}&=\frac{-K}{r}\log\left[\int\int\int_0^{+\infty}\eta_0^{b_1-1}\eta_1^{b_2-1}e^{\left[-\frac{\sum_{i=1}^n x_i}{\eta_0}-\sum_{i=1}^n\left(\frac{x_i}{\eta_1}\right)^{\beta}-a_1\eta_0-a_2\eta_1-r\eta_1\right]}\prod_{i=1}^n\left(\frac{1}{\eta_0}+\frac{\beta}{\eta_1}\left(\frac{x_i}{\eta_1}\right)^{\beta-1}\right)^{\delta_i}\text{d}\eta_0\text{d}\eta_1\text{d}\beta\right],\\
\hat{\beta}_{(L)}&=\frac{-K}{r}\log\left[\int\int\int_0^{+\infty}\eta_0^{b_1-1}\eta_1^{b_2-1}e^{\left[-\frac{\sum_{i=1}^n x_i}{\eta_0}-\sum_{i=1}^n\left(\frac{x_i}{\eta_1}\right)^{\beta}-a_1\eta_0-a_2\eta_1-r\beta\right]}\prod_{i=1}^n\left(\frac{1}{\eta_0}+\frac{\beta}{\eta_1}\left(\frac{x_i}{\eta_1}\right)^{\beta-1}\right)^{\delta_i}\text{d}\eta_0\text{d}\eta_1\text{d}\beta\right].
\end{align*}}
\end{table}
\noindent  The corresponding posterior risks are then
\[PR(\hat{\eta_0}_{(L)})=r(\hat{\eta}_{0(GQ)}-\hat{\eta}_{0(L)}),\]
\[PR(\hat{\eta}_{1(L)})=r(\hat{\eta}_{1(GQ)}-\hat{\eta}_{1(L)}),\]
\[PR(\hat{\beta}_{(L)})=r(\hat{\beta}_{(GQ)}-\hat{\beta}_{(L)}).\]
Since it is difficult to obtain closed form expressions of all these estimators, in the next section we will use the MCMC procedures to evaluate them.
\section{Simulation study}
In order to compare the performance of the  proposed Bayes estimators with the MLE estimators, we perform a Monte Carlo study assuming that $\eta_0=2, \eta_1=1 $ and  $\beta=2$ i.e. we consider the model $\mathcal{B}(2,1,2)$. Then, using $N=10000$ samples of the right censored model with different sizes $n=10, n=20$ and $n=30$. By choosing to  censor 10\%  respectively  20\% of date,  we obtain the following results.

\subsection{Likelihood estimation}
In the next tables we display the values of the estimators using  the EM algorithm for the $\mathcal{B}$ model when $10\%$ and $20\%$ of data are censored, where a Newton-Raphson algorithm is applied to the Weibull distribution and the direct likelihood maximization is applied to the exponential distribution.
\begin{table}[H]
\caption{The MLE of the parameters with quadratic error (in brackets)($10\%$).}
\begin{center}
\begin{tabular}{|c|c|c|} 
\hline
       $n$ & parameter & MLE \\
\hline
\multirow{3}{4cm}{\hspace{1.7cm}10} & $\eta_0$ & 1.8904 \;(0.0003)  \\
& $\eta_1$ &  06792\; (0.0710) \\
& $\beta$ & 1.9512 \; (0.1843) \\
\hline
\multirow{3}{4cm}{\hspace{1.7cm}20} & $\eta_0$ & 1.9158 \;(0.0085)  \\
& $\eta_1$ & 0.7913 \; (0.1031) \\
& $\beta$ & 1.9923 \; (0.0001) \\
\hline
\multirow{3}{4cm}{\hspace{1.7cm}30} & $\eta_0$ & 1.9985 \;(0.0002)  \\
& $\eta_1$ & 0.8181 \; (0.0011) \\
& $\beta$ & 1.9491 \; (0.0005) \\  
\hline     
\end{tabular}
\end{center}
\label{tab:ex}
\end{table}

\begin{table}[H]
\caption{The MLE of the parameters with quadratic error (in brackets)($20\%$).}
\begin{center}
\begin{tabular}{|c|c|c|} 
\hline
       $n$ & parameter & MLE \\
\hline
\multirow{3}{4cm}{\hspace{1.7cm}10} & $\eta_0$ & 1.9397 \;(0.0003)  \\
& $\eta_1$ & 0.6641 \; (0.0152) \\
& $\beta$ & 2.0005 \; (0.0001) \\
\hline
\multirow{3}{4cm}{\hspace{1.7cm}20} & $\eta_0$ & 1.9398 \;(0.0004)  \\
& $\eta_1$ & 0.9013 \; (0.0020) \\
& $\beta$ & 1.9491 \; (0.0001) \\
\hline
\multirow{3}{4cm}{\hspace{1.7cm}30} & $\eta_0$ & 1.8014 \;(0.0031)  \\
& $\eta_1$ & 0.9485 \; (0.0011) \\
& $\beta$ & 1.9611 \; (0.0030) \\  
\hline     
\end{tabular}
\end{center}
\label{tab:ex}
\end{table}

\subparagraph{Discussion:} For both censoring times, the estimated values of the parameters are close to the true values. Moreover, when 10\% of the data is censored, the smallest quadratic error corresponds to the largest $n$.

\subsection{Bayesian estimation}
The Bayesian estimators are obtained using the MCMC methods. For the choice of the hyperparameters designed from the equations given in Section 4. we consider the following prior informations.
For the shape parameter $\beta$, we assume that $[\beta_l,\beta_r]=[1,5]$, for the scale parameter $\eta_0$ of the exponential component we assume $[\eta_{0l},\eta_{0r}]=[1,300]$, and for the scale parameter $\eta_1$ of the Weibull component, we have $[\eta_{1l},\eta_{1r}]=[1,200]$. 

\begin{table}[H]
\caption{Bays estimators and PR (in brackets) under generalized quadratic loss function.}
\begin{center}
{\small
\begin{tabular}{*{9}{|c}|} 
\hline
       $n$ &\makecell{censoring\\percentage}& parameter & \multicolumn{6}{c|}{$\alpha$}\\
\cline{4-9}
 & & & -2&-1&-0.5&0.5&1&2\\ 
\hline
\multirow{2}{1cm}{\hspace{1cm}10} &\multirow{3}{1cm}{10$\%$} & \multirow{3}{2cm}{$\eta_0$}&2.0201&2.0619&2.0349&2.1191&2.1421&2.1600\\
& & &(0.0041)&(0.0072)&(0.0090) &(0.0082)& (0.0091)&(0.0091)\\
\cline{3-9}
& & \multirow{3}{2cm}{$\eta_1$}&1.1014&1.1014&1.1216&1.1425&1.1323&1.1338\\
& & &(0.0061)&(0.0711)&(0.0991)&(0.4005)&(0.2905)&(0.3136)\\
\cline{3-9}
& & \multirow{3}{2cm}{$\beta$}&1.9361&1.9132&1.8516&1.7331&1.7315&1.8405\\
& & &(0.0194)&(0.0621)&(0.8221)&(0.0872)&(0.0881)&(0.1416)\\
\cline{2-9}
&\multirow{3}{1cm}{20$\%$} & \multirow{3}{2cm}{$\eta_0$}&2.0717&2.0191&2.0301&2.0333&2.0509&2.0991\\
& & &(0.0049)& (0.0051)&(0.0059) &(0.0071)& (0.0071)&(0.0079)\\
\cline{3-9}
& & \multirow{3}{2cm}{$\eta_1$}&1.1609&1.0861&1.2615&1.5441&1.7822&1.9001\\
& & &(0.0612)&(0.0914)&(0.1009)&(0.0991)&(0.2923)&(0.4105)\\
\cline{3-9}
& & \multirow{3}{2cm}{$\beta$}&1.8306&1.8031&1.8094&1.8910&1.9700&2.0010 \\
& & &(0.0405)&(0.0511)&(0.0538) &(0.0711)& (0.0811)&(0.1009\\
\hline
\multirow{2}{1cm}{\hspace{1cm}20} &\multirow{3}{1cm}{10$\%$} & \multirow{3}{2cm}{$\eta_0$}&2.0001&2.0002&2.0001&1.9901&1.9491&1.5800\\
& & &(0.0001)&(0.0015)&(0.0018) &(0.0019)& (0.0049)&(0.0080)\\
\cline{3-9}
& & \multirow{3}{2cm}{$\eta_1$}&0.6615&0.6703&0.5301&0.5261&0.6401&0.7005 \\
& & &(0.0003)&(0.0021)&(0.0048)&(0.0014)&(0.0030)&(0.0003)\\
\cline{3-9}
& & \multirow{3}{2cm}{$\beta$}&1.9013&1.9005&1.9305&1.8001&1.6609&1.6712 \\
& & &(0.0001)&(0.0029)&(0.0005)&(0.0007)&(0.0039)&(0.0059)\\
\cline{2-9}
&\multirow{3}{1cm}{20$\%$} & \multirow{3}{2cm}{$\eta_0$} &2.0201&2.0482&2.0561&2.0531&2.0677&2.0823\\
& & &(0.0041)& (0.0070)&(0.0073) &(0.0061)& (0.0069)&(0.0072)\\
\cline{3-9}
& & \multirow{3}{2cm}{$\eta_1$}&0.9511&0.8552&0.8325&0.8512&0.8133&0.9205 \\
& & &(0.0051)&(0.0603)&(0.0991)&(0.1102)&(0.1512)&(0.2243)\\
\cline{3-9}
& & \multirow{3}{2cm}{$\beta$}&1.9431&1.9005&1.8522&1.8914&1.9233&1.9705 \\
& & &(0.0254)&(0.0491)&(0.0605)&(0.1231)&(0.2215)&(0.03105)\\
\hline
\multirow{2}{1cm}{\hspace{1cm}30} &\multirow{3}{1cm}{10$\%$} & \multirow{3}{2cm}{$\eta_0$}&2.1921&2.1905&2.1901&2.1883&2.1879&2.1863\\
& & &(0.0013)&(0.0014)&(0.0015) &(0.0017)& (0.0019)&(0.0021)\\
\cline{3-9}
& & \multirow{3}{2cm}{$\eta_1$}&1.1205&1.1201&1.1182&1.1173&1.1145&1.1129 \\
& & &(0.0004)&(0.0004)&(0.0005)&(0.0007)&(0.0029)&(0.0031)\\
\cline{3-9}
& & \multirow{3}{2cm}{$\beta$}&1.9421&1.9433&1.9441&1.9705&1.9733&1.9802 \\
& & &(0.0006)&(0.0006)&(0.0006)&(0.0007)&(0.0013)&(0.0023)\\
\cline{2-9}
&\multirow{3}{1cm}{20$\%$} & \multirow{3}{2cm}{$\eta_0$} &2.1031&2.0894&2.0972&2.0345&2.0372&2.0382\\
& & &(0.0013)& (0.0016)&(0.0014) &(0.0016)& (0.0017)&(0.0018)\\
\cline{3-9}
& & \multirow{3}{2cm}{$\eta_1$}&0.8305&0.8313&0.7805&0.7134&0.7235&0.8302 \\
& & &(0.0002)&(0.0003)&(0.0006)&(0.0007)&(0.0018)&(0.0034)\\
\cline{3-9}
& & \multirow{3}{2cm}{$\beta$}&1.9909&1.9733&1.9542&1.9521&1.9506&1.9503 \\
& & &(0.0005)&(0.0007)&(0.0008)&(0.0019)&(0.0032)&(0.0041)\\
\hline
\end{tabular}
}
\end{center}
\label{tab:ex}
\end{table}

\begin{table}[H]
\caption{Bays estimators and PR (in brackets) under the entropy loss function.}
\begin{center}
{\small
\begin{tabular}{*{9}{|c}|} 
\hline
       $n$ &\makecell{censoring\\percentage}& parameter & \multicolumn{6}{c|}{P} \\
\cline{4-9}
 & & & -2&-1&-0.5&0.5&1&2\\ 
\hline
\multirow{2}{1cm}{\hspace{1cm}10} &\multirow{3}{1cm}{10$\%$} & \multirow{3}{2cm}{$\eta_0$}&2.1093&2.1098&2.1046&2.0985&2.0941&2.0920\\
& & &(0.0020)&(0.0007)&(0.0081) &(0.0007)& (0.0019)&(0.0143)\\
\cline{3-9}
& & \multirow{3}{2cm}{$\eta_1$}&0.8113&0.7864&0.7182&0.6191&0.5914&0.5132 \\
& & &(0.1105)&(0.0258)&(0.0914)&(0.0601)&(0.3914)&(0.0181)\\
\cline{3-9}
& & \multirow{3}{2cm}{$\beta$}&1.9633&1.9103&1.8901&1.6515&1.7314&1.6105 \\
& & &(0.0532)&(0.0051)&(0.0301)&(0.0068)&(0.0313)&(0.1104)\\
\cline{2-9}
&\multirow{3}{1cm}{20$\%$} & \multirow{3}{2cm}{$\eta_0$} &2.1001&2.0931&2.1032&2.1013&2.0909&2.0891\\
& & &(0.0051)& (0.0009)&(0.0066) &(0.0006)& (0.0040)&(0.0021)\\
\cline{3-9}
& & \multirow{3}{2cm}{$\eta_1$}&0.8024&0.8005&0.7832&0.7214&0.6745&0.6691 \\
& & &(0.0532)&(0.0213)&(0.1154)&(0.1053)&(0.0713)&(0.0714)\\
\cline{3-9}
& & \multirow{3}{2cm}{$\beta$}&1.9501&1.8295&1.7917&1.7431&1.6913&1.6565 \\
& & &(0.0609)&(0.0071)&(0.0614)&(0.0615)&(0.0061)&(0.0324)\\
\hline

\multirow{2}{1cm}{\hspace{1cm}20} &\multirow{3}{1cm}{10$\%$} & \multirow{3}{2cm}{$\eta_0$}&2.1214&2.1032&2.1029&2.0963&2.0815&2.0803\\
& & &(0.0059)&(0.0061)&(0.0003) &(0.0070)& (0.0028)&(0.0203)\\
\cline{3-9}
& & \multirow{3}{2cm}{$\eta_1$}&1.0995&1.1015&1.1005&1.1001&1.0993&1.0957 \\
& & &(0.1414)&(0.0001)&(0.1405)&(0.0739)&(0.0729)&(0.1223)\\
\cline{3-9}
& & \multirow{3}{2cm}{$\beta$}&2.1818&2.1809&2.1793&2.1774&2.1751&2.1731 \\
& & &(0.0711)&(0.0005)&(0.0714)&(0.0089)&(0.2914)&(0.1095)\\
\cline{2-9}
&\multirow{3}{1cm}{20$\%$} & \multirow{3}{2cm}{$\eta_0$} &2.1781&2.1763&2.1751&2.1743&2.1731&2.1725\\
& & &(0.0008)& (0.0003)&(0.0004) &(0.0004)& (0.0004)&(0.0013)\\
\cline{3-9}
& & \multirow{3}{2cm}{$\eta_1$}&1.1535&1.1529&1.1719&1.1521&1.1509&1.1502 \\
& & &(0.0006)&(0.0006)&(0.0010)&(0.0007)&(0.0005)&(0.0032)\\
\cline{3-9}
& & \multirow{3}{2cm}{$\beta$}&2.1873&2.1859&2.1843&2.1839&2.1828&2.1819 \\
& & &(0.0013)&(0.0007)&(0.0022)&(0.0006)&(0.0006)&(0.0034)\\
\hline
\multirow{2}{1cm}{\hspace{1cm}30} &\multirow{3}{1cm}{10$\%$} & \multirow{3}{2cm}{$\eta_0$}&2.2011&2.2043&2.2032&2.2020&2.2018&2.2009\\
& & &(0.021)&(0.0003)&(0.0095) &(0.0004)& (0.0004)&(0.0039)\\
\cline{3-9}
& & \multirow{3}{2cm}{$\eta_1$}&1.1738&1.1729&1.1715&1.1709&1.1707&1.1702\\
& & &(0.0007)&(0.0003)&(0.0008)&(0.0006)&(0.0006)&(0.0020)\\
\cline{3-9}
& & \multirow{3}{2cm}{$\beta$}&2.1123&2.0832&2.0821&2.0819&2.0814&2.0809\\
& & &(0.0017)&(0.0003)&(0.0020)&(0.0004)&(0.0004)&(0.0020)\\
\cline{2-9}
&\multirow{3}{1cm}{20$\%$} & \multirow{3}{2cm}{$\eta_0$} &2.0829&2.0832&2.0821&2.0819&2.0814&2.0809\\
& & &(0.0007)& (0.0001)&(0.0020) &(0.0002)& (0.0006)&(0.0002)\\
\cline{3-9}
& & \multirow{3}{2cm}{$\eta_1$}&1.0734&1.0629&1.0583&1.0453&1.0423&1.0417 \\
& & &(0.0008)&(0.0005)&(0.0008)&(0.0006)&(0.0006)&(0.0001)\\
\cline{3-9}
& & \multirow{3}{2cm}{$\beta$}&1.9891&1.9877&1.9871&1.9868&1.9859&1.9843 \\
& & &(0.0008)&(0.0001)&(0.0003)&(0.0002)&(0.0005)&(0.0035)\\
\hline
\end{tabular}
}
\end{center}
\label{tab:ex}
\end{table}

\begin{table}[H]
\caption{Bays estimators and PR (in brackets) under the Linex loss function.}
\begin{center}
{\small
\begin{tabular}{*{9}{|c}|} 
\hline
       $n$ &\makecell{censoring\\percentage}& parameter & \multicolumn{6}{c|}{$r$} \\
\cline{4-9}
 & & &-2&-1&-0.5&0.5&1&2\\ 
\hline
\multirow{2}{1cm}{\hspace{1cm}10} &\multirow{3}{1cm}{10$\%$} & \multirow{3}{2cm}{$\eta_0$}&2.1333&2.1329&2.1284&2.1252&2.1134&2.1147\\
& & &(0.0182)&(0.0110)&(0.0009) &(0.0021)& (0.0091)&(0.0183)\\
\cline{3-9}
& & \multirow{3}{2cm}{$\eta_1$}&0.9784&0.9523&0.8613&0.7924&0.6813&0.5211 \\
& & &(0.3412)&(0.2656)&(0.0411)&(0.0420)&(0.3214)&(0.3929)\\
\cline{3-9}
& & \multirow{3}{2cm}{$\beta$}&2.0211&2.0220&1.9683&1.8967&1.8736&1.7969 \\
& & &(0.1523)&(1.1064)&(0.0209)&(0.0214)&(0.1241)&(0.1269)\\
\cline{2-9}
&\multirow{3}{1cm}{20$\%$} & \multirow{3}{2cm}{$\eta_0$} &2.0943&2.0937&2.0911&2.0901&2.0843&2.0823\\
& & &(0.0111)& (0.0084)&(0.0008) &(0.0009)& (0.2513)&(0.3409)\\
\cline{3-9}
& & \multirow{3}{2cm}{$\eta_1$}&0.9651&0.9117&0.7643&0.6853&0.6067&0.6095\\
& & &(0.3518)&(0.2614)&(0.0021)&(0.0527)&(0.2729)&(0.5114)\\
\cline{3-9}
& & \multirow{3}{2cm}{$\beta$}&1.9963&1.9695&1.9018&1.9001&1.8761&1.8569 \\
& & &(0.1563)&(0.0953)&(0.0214)&(0.0245)&(0.0305)&(0.1902)\\
\hline
\multirow{2}{1cm}{\hspace{1cm}20} &\multirow{3}{1cm}{10$\%$} & \multirow{3}{2cm}{$\eta_0$}&2.1426&2.0805&2.0793&2.0754&2.0721&2.0711\\
& & &(0.0113)&(0.0071)&(0.0008) &(0.0008)& (0.0021)&(0.0043)\\
\cline{3-9}
& & \multirow{3}{2cm}{$\eta_1$}&0.9374&0.8331&0.7852&0.6148&0.4747&0.4073 \\
& & &(0.4111)&(0.3119)&(0.0357)&(0.0434)&(0.3154)&(0.6189)\\
\cline{3-9}
& & \multirow{3}{2cm}{$\beta$}&2.1921&2.1920&2.1893&2.1884&2.1864&2.1859 \\
& & &(0.1616)&(0.0024)&(0.0003)&(0.0003)&(0.0020)&(0.0063)\\
\cline{2-9}
&\multirow{3}{1cm}{20$\%$} & \multirow{3}{2cm}{$\eta_0$} &2.1073&2.1051&2.1039&2.1021&2.1014&2.1009\\
& & &(0.0028)& (0.0023)&(0.0006) &(0.0007)& (0.0009)&(0.0018)\\
\cline{3-9}
& & \multirow{3}{2cm}{$\eta_1$}&1.1451&1.1443&1.1421&1.1419&1.1417&1.1405\\
& & &(0.0055)&(0.0052)&(0.0004)&(0.0005)&(0.0003)&(0.0002)\\
\cline{3-9}
& & \multirow{3}{2cm}{$\beta$}&2.0213&2.0205&2.0202&2.0183&2.0167&2.0166\\
& & &(0.0028)&(0.0032)&(0.0005)&(0.0006)&(0.0006)&(0.0020)\\
\hline
\multirow{2}{1cm}{\hspace{1cm}30} &\multirow{3}{1cm}{10$\%$} & \multirow{3}{2cm}{$\eta_0$}&2.1847&2.1833&2.1821&2.1819&2.1814&2.1809\\
& & &(0.0030)&(0.0350)&(0.0002) &(0.0003)& (0.0053)&(0.0061)\\
\cline{3-9}
& & \multirow{3}{2cm}{$\eta_1$}&1.1753&1.1733&1.1725&1.1.1801&1.1793&1.1722 \\
& & &(0.0052)&(0.0029)&(0.0005)&(0.0059)&(0.0029)&(0.0020)\\
\cline{3-9}
& & \multirow{3}{2cm}{$\beta$}&2.1853&2.1847&2.1835&2.1820&2.1825&2.1819 \\
& & &(0.0039)&(0.0027)&(0.0005)&(0.0006)&(0.0031)&(0.0037)\\
\cline{2-9}
&\multirow{3}{1cm}{20$\%$} & \multirow{3}{2cm}{$\eta_0$}&2.1653&2.1647&2.1635&2.1641&2.1629&2.1617\\
& & &(0.0021)& (0.0020)&(0.0004) &(0.0004)& (0.0016)&(0.0017)\\
\cline{3-9}
& & \multirow{3}{2cm}{$\eta_1$}& 0.9872&0.8822&0.7581&0.5979&0.4421&1.3771\\
& & &(0.0057)&(0.0023)&(0.0040)&(0.0004)&(0.0202)&(0.0075)\\
\cline{3-9}
& & \multirow{3}{2cm}{$\beta$}&2.2094&2.2087&2.2073&2.2059&2.2049&2.2038\\
& & &(0.0157)&(0.0023)&(0.0004)&(0.0014)&(0.0032)&(0.0055)\\
\hline
\end{tabular}
}
\end{center}
\label{tab:ex}
\end{table}

\subparagraph{Discussion:}
We note that the value $\alpha=-2$ gives us the best posterior risk associated with  generalised quadratic loss function. In the two censuring times, we also obtain the smallest risk when $n$ is large. Under the entropy loss function, we obtain the best posterior risk when $p=-1$ and $n=30$. Finally, under the Linex loss function, the case $r=-0.5$ provides the best results.
\\

\vspace{0.3cm}The following table illustrates the Bayesian estimator under the three loss functions.
\begin{table}[H]
\caption{Bays estimators and PR (in brackets) under the three loss function}
\begin{center}
{\small
\begin{tabular}{*{6}{|c}|} 
\hline
       $n$ &\makecell{censoring \\ percentage}& parameter &GQ$(\alpha=-2)$ & entropy $(p=-1)$&Linex $(r=-0.5)$ \\
\hline
\multirow{2}{1cm}{\hspace{1cm}10} &\multirow{3}{1cm}{10$\%$} & $\eta_0$&2.0201(0.0041) &2.1098 (0.0007)&2.1284(0.0009)\\
\cline{3-6}
& & $\eta_1$&1.1014(0.0061) &0.7864(0.0258)&0.8613(0.0411)\\
\cline{3-6}
& & $\beta$& 1.9361(0.0194)& 1.9103(0.0051)&1.9683(0.0209)\\
\cline{2-6}
&\multirow{3}{1cm}{20$\%$}& $\eta_0$&2.0931(0.0009) &2.0954(0.0008)&2.0911(0.0008) \\
\cline{3-6}
& & $\eta_1$&1.1609(0.0612) &0.8005(0.0213)&0.7643(0.0021) \\
\cline{3-6}
& & $\beta$&1.8306(0.0405) &1.8295(0.0071)&1.9018(0.0314)\\

\hline
\multirow{2}{1cm}{\hspace{1cm}20} &\multirow{3}{1cm}{10$\%$} & $\eta_0$&2.0001(0.0001) &2.1032(0.0061)&2.0793(0.0008)\\
\cline{3-6}
& & $\eta_1$& 0.6615(0.0003)&1.1015(0.0001)&0.7852(0.0375)\\
\cline{3-6}
& & $\beta$&1.9013(0.0001) &2.1809(0.0005)&2.1893(0.0003) \\
\cline{2-6}
&\multirow{3}{1cm}{20$\%$} & $\eta_0$ &2.0201(0.0041) &2.1763(0.0003)&2.1039(0.0006)\\
\cline{3-6}
& & $\eta_1$&0.9511(0.0051) & 1.1529(0.0006)&1.1421(0.0004)\\
\cline{3-6}
& & $\beta$&1.9381(0.0254) & 2.1859(0.0007)&2.0202(0.0005)\\
\hline
\multirow{2}{1cm}{\hspace{1cm}30} &\multirow{3}{1cm}{10$\%$} & $\eta_0$&2.1921(0.0013) &2.2043(0.0003)&2.1821(0.0002)\\
\cline{3-6}
& & $\eta_1$&1.1205(0.0004) &1.1729(0.0032)&1.1725(0.0005) \\
\cline{3-6}
& & $\beta$&1.9421(0.0006) & 2.1071(0.0003)&2.1835(0.0005)\\
\cline{2-6}
&\multirow{3}{1cm}{20$\%$} & $\eta_0$ &2.1031(0.0013) &2.0832(0.0001)&2.1635(0.0004)\\
\cline{3-6}
& & $\eta_1$&0.8305(0.0002) &1.0629(0.0005)&0.7581(0.0040) \\
\cline{3-6}
& & $\beta$ &1.9909(0.0005) &1.9877(0.0001)&2.2073(0.0004) \\
\hline
\end{tabular}
}
\end{center}
\label{tab:ex}
\end{table}
\subparagraph{Discussion:} We notice that the entropy loss function provides the best Bayesian estimator of the parameters among the three loss functions in regard to the posterior risk values. While clearly the other two (GQ, Linex) gave the same level of performance.
\subsection{The Survival and the Hazard functions}
Since the survival and the hazard functions are both depend on time, we consider here the interval of time $t=[1,50]$, then we obtain the following graphs representing the curves of the survival function with the real values of the parameters, and the estimated ones.

We consider the sampling size $n=30$, the MLE and the Bayesian estimators, we choose the estimator under the entropy loss function (the best estimator as we have seen above). in the both cases of censoring times  $10\%$ and $20\%$.  
\begin{figure}[H]
\includegraphics{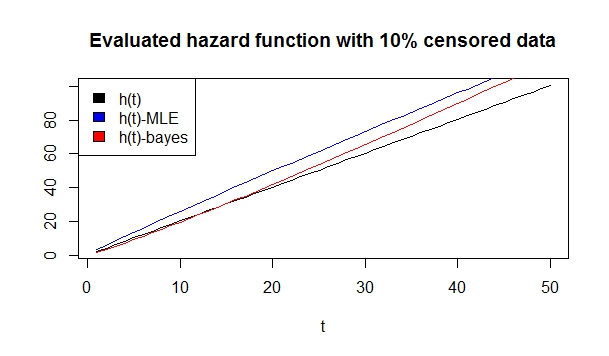}
\end{figure}
\begin{figure}[H]
\includegraphics{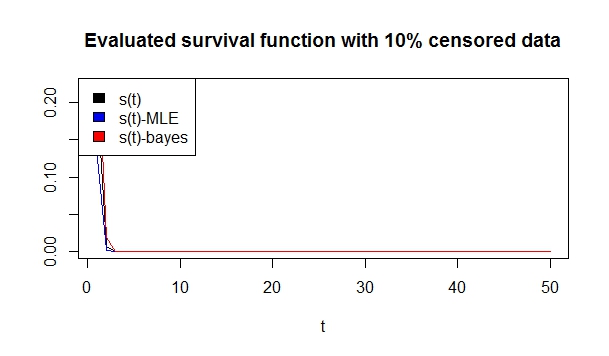}
\end{figure}
\begin{figure}[H]
\includegraphics{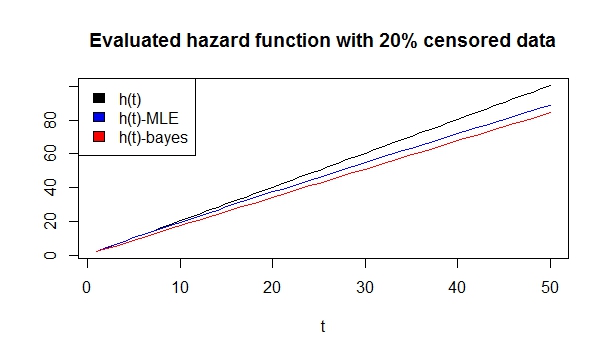}
\begin{figure}[H]
\includegraphics{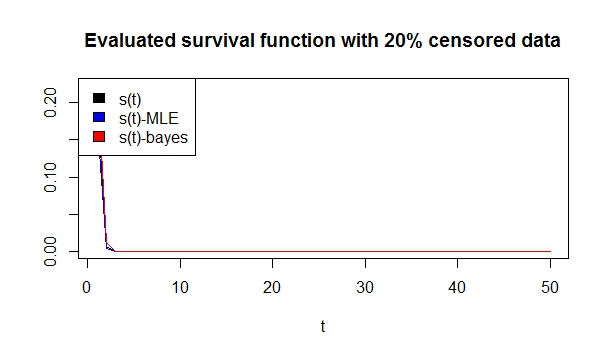}
\end{figure}
\end{figure}
\subparagraph{Discussion:}
We notice that in the case of the survival function, when $t<5$ the three functions are close to each other, while they do coincide when $t>5$, otherwise in the case of the hazard function, it is clearly that the Bayesian estimator performs better when $10\%$ of data are censored , in other words the survival function using the Bayesian estimates is closer to the survival function using the real values, and we notice the opposite case (the MLE estimator performs better) when $20\%$ of data are censured.

\subsection{Comparison of the estimators}
In this section, we compare the best Bayesian estimators obtained above with the
maximum likelihood estimators. For this, we propose to use the following two criteria: the Pitman closeness criterion
(see Pitman (1937), Fuller (1982) and Jozani (2012)) and the integrated mean square error (IMSE) defined as follows.

\subsubsection{Definition} 
An estimator  $\theta_1$ of a parameter $\theta$ dominates  another estimator  $\theta_2$ in the sense of Pitman closeness criterion if for all $\theta\in\Theta$
$$
P_{\theta}[|\theta_1-\theta|<|\theta_2-\theta|]>\frac{1}{2}.
$$
\subsubsection{Definition} Consider the estimators $\theta_i (i=1,\ldots,N)$ obtained with $N$ samples of model. The integrated mean square error is defined as
$$
\text{IMSE}=(\sum_{i=1}^N(\theta_i-\theta)^2 )/N.
$$

In the following tables, we present the values of the Pitman probabilities which allow us to compare the Bayesian estimators with the MLE under the three loss function when  $\alpha=-2, p=-1$ and  $r=-\frac{1}{2}$. 

\begin{table}[H]
\caption{Pitman comparison of the  estimators}
\begin{center}
{\small
\begin{tabular}{*{6}{|c}|} 
\hline
       $n$ &\makecell{censoring\\percentage}&parameter&GQ$(\alpha=-2)$ &entropy $(p=-1)$ b&Linex $(r=-\frac{1}{2})$ \\
\hline
\multirow{2}{1cm}{\hspace{1cm}10} &\multirow{3}{1cm}{10$\%$} & $\eta_0$&0.735 &0.719& 0.682\\
\cline{3-6}
& & $\eta_1$&0.352 &0.205&0.226\\
\cline{3-6}
& & $\beta$&0.682 &0.575 &0.594\\
\cline{2-6}
&\multirow{3}{1cm}{20$\%$}& $\eta_0$&0.546 &0.522&0.519 \\
\cline{3-6}
& & $\eta_1$&0.366&0.312&0.299\\
\cline{3-6}
& & $\beta$&0.557&0.519&0.533 \\
\hline
\multirow{2}{1cm}{\hspace{1cm}20} &\multirow{3}{1cm}{10$\%$} & $\eta_0$&0.134&0.148&0.175\\
\cline{3-6}
& & $\eta_1$&0.566&0.593&0.629\\
\cline{3-6}
& & $\beta$&0.335&0.394&0.341 \\
\cline{2-6}
&\multirow{3}{1cm}{20$\%$} & $\eta_0$ &0.288& 0.275&2.242\\
\cline{3-6}
& & $\eta_1$& 0.501&0.516&0.513\\
\cline{3-6}
& & $\beta$&0.205&0.201&0.197\\
\hline
\multirow{2}{1cm}{\hspace{1cm}30} &\multirow{3}{1cm}{10$\%$} & $\eta_0$&0.349&0.318&0.299\\
\cline{3-6}
& & $\eta_1$&0.561&0.601&0.621 \\
\cline{3-6}
& & $\beta$&0.127&0.122&0.119\\
\cline{2-6}
&\multirow{3}{1cm}{20$\%$} & $\eta_0$ &0.115&0.118&0.202\\
\cline{3-6}
& & $\eta_1$& 0.231&0.209&0.285\\
\cline{3-6}
& & $\beta$ &0.243&0.219&0.208\\
\hline
\end{tabular}
}
\end{center}
\label{tab:ex}
\end{table}
\subparagraph{Discussion:} When $n$ is small, the Bayesian estimators of $\eta_0$ and $\beta$ are better than the MLE estimators. Furthermore, we  note that the generalised quadratic loss function provides the best values, however the MLE estimator of $\eta_1$ is better than the Bayesian one. When $n$ is large, the MLE estimators of the three parameters are better than the Bayesian estimators.

\noindent In the next table we present the values of the integrated mean square error of the estimators under the three loss function and the maximum likelihood estimator. 
\begin{table}[H]
\caption{The IMSE of the estimators}
\begin{center}
{\small
\begin{tabular}{*{7}{|c}|} 
\hline
       $n$ &\makecell{censoring\\percentage}& parameter&MLE &GQ$(\alpha=-2)$&entropy $(p=-1)$&Linex $(r=-0.5)$ \\
\hline
\multirow{2}{1cm}{\hspace{1cm}10} &\multirow{3}{1cm}{10$\%$} & $\eta_0$ & 0.1481& 0.0061& 0.0106 & 0.0081\\
\cline{3-7}
& & $\eta_1$&0.0053 &0.2914 &0.2143&0.2205\\
\cline{3-7}
& & $\beta$&0.1934&0.0714 &0.0721&0.0034\\
\cline{2-7}
&\multirow{3}{1cm}{20$\%$}& $\eta_0$&0.0713 &0.0007&0.0031&0.0033 \\
\cline{3-7}
& & $\eta_1$&0.0529&0.2305&0.2115&0.2101\\
\cline{3-7}
& & $\beta$&0.0105&0.0715&0.0063&0.0042 \\
\hline
\multirow{2}{1cm}{\hspace{1cm}20} &\multirow{3}{1cm}{10$\%$} & $\eta_0$&0.0613&0.0043&0.0073&0.0031\\
\cline{3-7}
& & $\eta_1$&0.0359&0.2093&0.2063&0.2215 \\
\cline{3-7}
& & $\beta$& 0.1048&0.0011&0.0053&0.0043\\
\cline{2-7}
&\multirow{3}{1cm}{20$\%$} & $\eta_0$ &0.0801&0.2301&0.2297&0.2543\\
\cline{3-7}
& & $\eta_1$&0.0563&0.0984&0.1001&0.0975\\
\cline{3-7}
& & $\beta$&0.0421&0.0405&0.0463&0.0441 \\
\hline
\multirow{2}{1cm}{\hspace{1cm}30} &\multirow{3}{1cm}{10$\%$} & $\eta_0$&0.1622&0.0063&0.0501&0.0479\\
\cline{3-7}
& & $\eta_1$& 0.0393&0.0152&0.0322&0.0310\\
\cline{3-7}
& & $\beta$& 0.1682&0.0215&0.0329&0.0308\\
\cline{2-7}
&\multirow{3}{1cm}{20$\%$} & $\eta_0$ &0.0751&0.0511&0.0521&0.0522\\
\cline{3-7}
& & $\eta_1$&0.4242&0.3952&0.4102&0.4088 \\
\cline{3-7}
& & $\beta$ & 0.2905 & 0.2143 & 0.2184 & 0.2123\\
\hline
\end{tabular}
}
\end{center}
\label{tab:ex}
\end{table}
\subparagraph{Discussion:} 
When $n$ is small, the Bayesian estimators of $\eta_0$ and $\beta$ provides the smallest IMSE compared to the MLE estimators. But in for $\eta_1$ the MLE estimator preforms better than the Bayesian one. When $n$ is large, all the Bayesian estimators are better than the MLE estimator, and we can notice that the generalised quadratic loss function provide the best values of the IMSE.

\section{Conclusion}
In this study we considered a simple competing risk model based on Weibull and exponential failures, We used classical and Bayesian estimation methods to estimate the unknown parameters where we used  the EM algorithm since no closed form of the MLE estimators can be obtained.  The results were obtained using simulated data sets of size 10, 20 and 30. We obtained the Bayesian estimators under the generalized quadratic, entropy and Linex loss functions.Then we used the Monte-Carlo simulation technique to determine which loss function has the smallest posterior risks. These selected Bayesian estimators are compared with the maximum likelihood estimators of the unknown parameters using Pitman's closeness criterion and the integrated mean square error. as  future prospect,  a mixture of the loss functions used in this paper might yield an optimal estimation.

\section*{References}

\bibliography{mybibfile}

\end{document}